\documentclass[12pt,a4paper]{article}
\usepackage{amsmath}
\usepackage{amssymb}
\usepackage{mathrsfs}
\usepackage{indentfirst}
\usepackage{geometry}
\usepackage{amscd}
\usepackage{amsthm}
\usepackage{latexsym}
\usepackage{amsfonts}
\usepackage{epsfig}
\usepackage{color}
\usepackage{graphicx}
\usepackage{array}
\begin{document}

\renewcommand{\a}{\alpha}
\newcommand{\D}{\Delta}
\newcommand{\ddt}{\frac{d}{dt}}
\numberwithin{equation}{section}
\newcommand{\e}{\epsilon}
\newcommand{\eps}{\varepsilon}

\newtheorem{theorem}{Theorem}[section]
\newtheorem{proposition}[theorem]{Proposition}
\newtheorem{corollary}[theorem]{Corollary}
\newtheorem{lemma}[theorem]{Lemma}
\newtheorem{definition}[theorem]{Definition}
\newtheorem{example}[theorem]{Example}
\newtheorem{remark}[theorem]{Remark}
\newcommand{\dis}{\displaystyle}
\newtheorem{examples}[theorem]{Examples}
\renewcommand{\thefootnote}{\fnsymbol{footnote}}

\title{\bf Inverse scattering for the nonlinear magnetic Schr\"{o}dinger equation\footnote{This work was supported by National Natural Science Foundation of China (61671009, 12171178).}}
\author{Lei Wei $^{\star\,1}$\,and\ Hua Huang$^{2}$\,\\
{\small {\it $^{1}$School of Information and Artificial Intelligence, Anhui Agricultural University,}}\\
{\small {\it Hefei, {\rm 230036,} P.R.China.}}\\
{\small {\it $^{2}$School of Information Engineering, Wuhan College, Wuhan, P.R.China.}}\\
{\small {\it $^{\star}$Corresponding author(s). Email:  weileimath@hust.edu.cn}}\\
{\small {\it Contributing authors. Email:huahuang@whxy.edu.cn}}}
\date{}
\maketitle

\noindent{\large\bf Abstract}\\
  
In this paper, we focus on the inverse scattering problem for the nonlinear Schr\"odinger equation with magnetic potentials. Specifically, we investigate whether the scattering operator associated with the nonlinear Schr\"odinger equation can uniquely determine the magnetic potential. Our main goal is to establish the uniqueness result for the magnetic potential based on the scattering data obtained from the scattering operator. \\



\noindent{\bf Authors’ Contributions:} All authors contributed equally to this work. All authors read and approved the final manuscript.\\

%\noindent{\bf Funding Statement:} This work was supported by National Science Foundation of China (NSFC) (61671009,12171178).\\

\noindent{\bf Conflicts of Interest:} The authors declare that they have no conflicts of interest to report regarding the present study.

\section{Introduction}
We consider the nonlinear magnetic Schr\"odinger equation 
\begin{equation}\label{equation1}
\left\{
\begin{split}
&(i\partial_{t}+H+V)u=\arrowvert u\arrowvert^{p-1}u, \ x\in \mathbb{R}^{n},\ t\in\mathbb{R},\\
&u(0)=\phi,
\end{split}
\right.
\end{equation}
where $H:=(\nabla+iA)^{2}=\Delta+iA\cdot\nabla+i\nabla\cdot A-\arrowvert A\arrowvert^{2}$,  $A(x)=(A_{1}(x),A_{2}(x),...,A_{n}(x))$ and $V(x)$ are magnetic vector and electric scalar potentials of the field.

In all the paper, our assumptions will be the following: \\
(\uppercase\expandafter{\romannumeral1}) Assume that the magnetic Schr\"{o}dinger operator $H$ has the following decay estimates,
\begin{flalign}\label{assum1}
\|e^{itH}\phi\|_{L^\infty}\leqslant C|t|^{-n/2}\|\phi\|_{L^1}.
\end{flalign}
(\uppercase\expandafter{\romannumeral2}) Assume that the magnetic potential is exponentially decaying, i.e.,
\begin{flalign}\label{assum2}
A\in e^{-\gamma_0\langle x\rangle}W^{1,\infty}(\mathbb{R}^n,\,\mathbb{R}^n),
\end{flalign}
where $\gamma_0>0,\langle x\rangle=(1+|x|^2)^{\frac{1}{2}}$.\\
(\uppercase\expandafter{\romannumeral3})
Assume that for $j\in{1, ... ,n}$, $A_{j}(x)$ is real valued, $C^{\infty}$ on $\mathbb{R}^{n}$. If $B = (B_{jk})$ with $B_{jk}= \partial_{j}A_{k}- \partial_{k}A_{j}$, then there exists $\varepsilon>0$ such that
\begin{flalign}\label{assum3}
&\arrowvert\partial^{\alpha}B(x)\arrowvert\leqslant C_{\alpha}(1+\arrowvert x\arrowvert)^{-n-\varepsilon},\nonumber\\
&\arrowvert\partial^{\alpha}A(x)\arrowvert\leqslant C_{\alpha},\ \forall\alpha\geqslant1,\ \forall x\in\mathbb{R}^n.
\end{flalign}
In addition for $V(x)$, we assume that $\arrowvert\partial^{\alpha}V(x)\arrowvert\leqslant C_{\alpha},\ \forall\alpha\geqslant2$ and there exists $m>0$ such that $V(x)\geqslant m,\ \forall x\in\mathbb{R}^n$. 

The goal of this paper is to give a method to uniquely reconstruct the potentials $A(x)$ and $V(x)$ from the scattering operator associated with \eqref{equation1} and the corresponding unperturbed linear Schr\"odinger equation, 
\begin{equation}\label{equation2}
\left\{
\begin{split}
&(i\partial_{t}+H_{0})u=0, \ x\in \mathbb{R}^{n},\ t\in\mathbb{R},\\
&u(0)=\phi,
\end{split}
\right.
\end{equation}
where $H_{0}=\Delta$.

Let's recall the well-posedness of solutions for the magnetic Schr\"{o}dinger equations. The well-posedness of weak solutions of the nonlinear Schr\"odinger equation was obtained in the case of $A=0$ and $V=0$ in \cite{K1,K2,T1}, and in the case of $A\neq0$ and $V\neq0$ in  \cite{B,Na}. For the local well-posedness of strong solutions of the same equation, it was acquired in \cite{K2,T2} when $A=0$, and in \cite{NS,WD} when $A\neq0$. 

For the scattering, the nonlinear Schr\"odinger equation with $A=0$ and $V=0$, has been extensively studied by many mathematicians and physicists almost since its advent. Struss \cite{S,S1} studied the scattering theory and applyed it to Schr\"odinger and Klein-Gordon equations with various kinds of nonlinear terms. Bourgain \cite{Bo} established the scattering in the energy space $H^{1}(\mathbb{R}^{3})$ for defocusing nonlinear Schr\"odinger equation. In the well know case $A=0$ and $V\neq0$ for linear Schr\"odinger equation, the term \textit{short range} physically means potential $V(x)$ decays fast enough to ensure that a quantum scattering system has locality in a large scale. See Agmon \cite{Ag} and H\"{o}rmander \cite{H} proved the spectral and scattering theory when $V(x)$ satisfies a short range condition, i.e. $V(x)$ decays at infinity like $\sim \arrowvert x\arrowvert^{-1-\varepsilon},\varepsilon>0$. In the recent survey, Schlag \cite{Sc2} established a structure formula for the intertwining wave operators in $\mathbb{R}^{3}$, which as applied to spectral theory. For the magnetic Schr\"{o}dinger equation, there are not many results here, what is known is Demuth and Ouhabaz \cite{DO} proved the existence and completeness of the wave operators $W_{\pm}(A(b),-\Delta)$ with $A(b)=-(\nabla-ib(x))^{2}$, $x\in\mathbb{R}^{n}$. 

For the inverse scattering problem, Enss and Weder \cite{EW1, EW2} establiesd a new time-dependent method for the linear Schr\"odinger equation with $A=0$ but $V\neq0$. For scattering in the presence of a constant electric field see Weder \cite{W2}. A key idea of the method of \cite{EW1, EW2} and \cite{W2} is to consider high velocity scattering states. Then Weder \cite{W, W1} proved the inverse scattering for nonlinear free Schr\"odinger equation $A=0$ and $V=0$. For the inverse scattering problems of the magnetic Schr\"{o}dinger operator, the magnetic fields and electric potentials were uniquely determined in \cite{ER} which considered potentials with many exponentially decaying derivatives, and this results were improved in \cite{PSU} for exponentially decaying potentials. In this paper, with different technology from \cite{ER,PSU}, our main method of this paper is inspired by \cite{W, W1}. That is to say, we need to modulate the amplitude of the scattering state which allows us to uniquely reconstruct the potential considering scattering states that have small amplitude and high velocity. 

\begin{remark}\label{remark1}
For dispersive estimate \eqref{assum1} of the magnetic nonlinear Schr\"odinger equations, it's an open problem until now. But we need the decay of $t$ to insure the existence of nonlinear scattering problem (See Theorem \ref{theorem1.2}), then we consider the inverse problem (See Theorem \ref{theorem1.3} and Corollary \ref{corollary1.4}).
\end{remark}

Before consider the inverse scattering problems, we have to state some general nonlinear scattering results. Similar to Strauss \cite{S,S1}, we have the existence of the scattering operator for \eqref{equation1}. Denote $H_\delta=\{\phi\in H^1, ||\phi||_{H^1}<\delta\}$.

\begin{theorem}\label{theorem1.2}
For Schr\"odinger equation \eqref{equation1} with $V=0$ in $n\geqslant3$, we assume that Schr\"odinger operator $H$ and potential $A(x)$ satisfy the hypothesis (\uppercase\expandafter{\romannumeral1}) and (\uppercase\expandafter{\romannumeral2}). If  $\phi_-\in H_\delta$ for some $\delta>0$,  then there is a unique solution of the integral equation
\begin{flalign}\label{1.2.1}
u(t)=e^{itH}\phi_-+\int^t_{-\infty}e^{i(t-\tau)H}\arrowvert u\arrowvert^{p-1}u(\tau)d\tau
\end{flalign}
such that
$$u\in C(\mathbb{R},\,H^1)\cap L^{r}(\mathbb{R},\,W^{k,1+p})$$
and
\begin{flalign}\label{1.2.2}
\lim_{t\rightarrow -\infty}||u(t,\cdot)-e^{itH}\phi_-||_{H_1}=0.
\end{flalign}
Furthermore, there exists a unique $\phi_+\in H_\delta$ such that
\begin{flalign}\label{1.2.3}
\lim_{t\rightarrow +\infty}||u(t,\cdot)-e^{itH}\phi_+||_{H_1}=0.
\end{flalign}
In addition,
\begin{flalign}\label{1.2.4}
||u(t,x)||_{L^2} \ \ and \ \ E(u)=\int  [\frac{1}{2}|\sqrt{H}u(t,\cdot)|^2+\arrowvert u(t,\cdot)\arrowvert^{p+1}]dx
\end{flalign}
are finite and independent of $t$, and
\begin{flalign}\label{1.2.5}
u(t)=e^{itH}\phi_+-\int^{+\infty}_te^{i(t-\tau)H}\arrowvert u\arrowvert^{p-1}u(\tau)d\tau.
\end{flalign}
\end{theorem}

From Theorem \ref{theorem1.2}, the scattering operator is well-defined on some neighborhood of zero in $H^1$ for the equation \eqref{equation1}. For a proof, see \cite{HT,S1}. By \eqref{1.2.1} and \eqref{1.2.5} we have
\begin{flalign}\label{1-1}
\phi_+=\phi_-+\int_{-\infty}^{+\infty}e^{-i\tau H}\arrowvert u\arrowvert^{p-1}u(\tau)d\tau.
\end{flalign}
Then scattering operator $S_{A}$, defined by $S_{A}:\ \phi_-\rightarrow \phi_+$ is a homeomorphism from $H_\delta$ onto $H_\delta$.

Denote by $S_{L}$ the linear scattering operator that compares solutions to linear Schr\"odinger equation \eqref{equation1} without nonlinear term with solutions to linear Schr\"odinger equation \eqref{equation2}:
\begin{flalign}\label{SL}
S_{L}=W^{*}_{+}W_{-}.
\end{flalign}
In the following Theorem \ref{theorem1.3}, we reconstruct the linear scattering operator $S_{L}$ from the small amplitude limit of the nonlinear scattering operator $S_{A}$. 
\begin{theorem}\label{theorem1.3}
For Schr\"odinger equation \eqref{equation1} with $V=0$ in $n\geqslant3$, we assume that Schr\"odinger operator $H$ and potential $A(x)$ satisfy the hypothesis (\uppercase\expandafter{\romannumeral1}) and (\uppercase\expandafter{\romannumeral2}). Then for any $\phi,\ \psi \in H^{1}$,
\begin{flalign}\label{1.3.1}
\lim_{\epsilon\downarrow 0}\frac{1}{\epsilon}(S_{A}\epsilon\phi,\psi)=(S_L\phi,\psi).
\end{flalign}
\end{theorem}

From above we can uniquely reconstruct $S_L$ from the small amplitude limit of $S$.
Since we can rewrite the operator $H$ as the form $H=\sum_{j=1}^n(D_j+A_j)^2-\arrowvert A\arrowvert^{2}$,
then we apply a inverse scattering result of the magnetic Schr\"{o}dinger operator from \cite{PSU},
which stats that if the scattering matrices for two sets of exponentially decaying coefficients coincide at a fixed energy,
then the magnetic fields and electric potentials have to be the same.
In this way we have the following Corollary:

\begin{corollary}\label{corollary1.4}
From Theorem \ref{theorem1.3}, assume that (\uppercase\expandafter{\romannumeral1}) and (\uppercase\expandafter{\romannumeral2}) are satisfied, then the scattering operator $S_{A}$ determines uniquely the the magnetic potential $A(x)$.
\end{corollary}

\begin{theorem}\label{theorem1.5}
For Schr\"odinger equation \eqref{equation1} in $n\geqslant2$, we assume that potentials $A(x)$ and $V(x)$ satisfy the hypothesis (\uppercase\expandafter{\romannumeral3}). Then we have for any $\phi,\ \psi \in \Sigma$,
\begin{flalign}\label{1.5.1}
\lim_{\epsilon\downarrow 0}\frac{1}{\epsilon}(S\epsilon\phi,\psi)=(S_L\phi,\psi),
\end{flalign}
where $\Sigma=\{\phi\in\mathscr{S'}(\mathbb{R}^{n}),\ (-\Delta)^{\frac{s}{2}}\phi\in L^{2}(\mathbb{R}^{n}),\ \arrowvert x\arrowvert^{s}\phi\in L^{2}(\mathbb{R}^{n}),\ |x|^{s_{2}}(-\Delta)^{\frac{s_{1}}{2}}\phi\in L^{2}(\mathbb{R}^{n})\}$, $s>\frac{n}{2}$ and $s=s_{1}+s_{2}$, $0<s_{1},s_{2}<s$.
\end{theorem}

\begin{corollary}\label{corollary1.6}
From Theorem \ref{theorem1.5}, assume that (\uppercase\expandafter{\romannumeral3}) is satisfied, then the scattering operator $S$ determines uniquely the potentials $A(x)$ and $V(x)$.
\end{corollary}

Our paper is organized as follows. In Section 2, we mainly focus on proving Theorem \ref{theorem1.3}. In Section 3, we concentrate on proving Theorem \ref{theorem1.5}. The method used in the proofs of both Theorem \ref{theorem1.3} and Theorem \ref{theorem1.5} are adapted from \cite{W,W1}, but the original method is not directly applicable to the magnetic Schrödinger operator. Therefore, we have made some modifications and adjustments during the proof.\\

\addcontentsline{toc}{section}{Notation}
\hspace{-5mm}\textbf{Notation}:
\begin{itemize}
\setlength{\baselineskip}{1.5\baselineskip}
	
\item $\mathscr{F}u(x)=\hat{u}=(2\pi)^{-\frac{n}{2}}\int_{\mathbb{R}^{n}}e^{-ix\cdot \xi}u(\xi)d\xi$.

\item
$(\mathscr{F}^{-1}u)(x)=\check{u}=(2\pi)^{-\frac{n}{2}}\int_{\mathbb{R}^{n}}e^{ix\cdot \xi}u(\xi)d\xi$.

\item
$H^{s}(\mathbb{R}^{n})=\{u\in L^{2}(\mathbb{R}^{n}):(1+\arrowvert\xi\arrowvert^{2})^{\frac{s}{2}}\hat{u}\in L^{2}(\mathbb{R}^{n})\}$.

\item
$W^{k,p}(\mathbb{R}^{n})=\{u\in L^{p}(\mathbb{R}^{n}):(1+\arrowvert\xi\arrowvert^{2})^{\frac{s}{2}}\hat{u}\in H^{k}(\mathbb{R}^{n})\}$.

\end{itemize}	

\section{Proof of Theorem \ref{theorem1.3}}
Before prove, we state some results from \cite{JSS,S,S1} which I need after. Denote $Z=L^r(\mathbb{R},\,L^{1+p})\cap L^{\infty}(\mathbb{R},\,L^{1+p})$, where $(r,1+p)\in\wedge$, and the norm of $Z$ is given by
\begin{flalign*}
\|u\|_Z=\|u\|_{L^r(\mathbb{R},\,L^{1+p})}+\|u\|_{L^{\infty}(\mathbb{R},\,L^{1+p})}.
\end{flalign*}
Denote
\begin{flalign*}
\mathscr{P}_su(t):=\int_s^te^{i(t-\tau)H}\arrowvert u\arrowvert^{p-1}u(\tau)d\tau,
\end{flalign*}
where $u(\tau)\in Z$.

From the formal representation the nonlinear term, we have for $\phi,\psi\in L^{1+p}$
\begin{flalign}\label{2.1}
\|\arrowvert\phi\arrowvert^{p-1}\phi-\arrowvert\psi\arrowvert^{p-1}\psi\|_{L^{1+p}}\leq C(\|\phi\|_{L^{1+p}}+\|\psi\|_{L^{1+p}})^{p-1}\|\phi-\psi\|_{L^{1+p}},
\end{flalign}
then using the weak Young's inequality, it follows that
\begin{flalign}\label{2.2}
\|\mathscr{P}_s\phi-\mathscr{P}_s\psi\|_{L^{1+p}}\leq C(\|\phi\|_{L^{1+p}}+\|\psi\|_{L^{1+p}})^{p-1}\|\phi-\psi\|_{L^{1+p}}.
\end{flalign}
Moreover, note the fact $\mathscr{P}_su(t)\in C(R,L^{1+p})$. Let's denote $Z(\delta)=\{u\in Z, \|u\|_Z\leq \delta\}$. Since $p>1$, by \eqref{2.2} and assumption \eqref{assum1}, we can obtain that $\mathscr{P}_s$ is a contraction on $Z(\delta)$ if $\delta$ is small enough. In fact for the specific proof, see \cite{S}. 

Meanwhile, it follows from the estimate \eqref{assum1} that
\begin{flalign}\label{2.3-1}
\|e^{itH}\phi\|_{L^r(\mathbb{R},\,L^{1+p})} \leq C\|\phi\|_{H^1}.
\end{flalign}
Note that $H^{1}\subset L^{1+p}$,  $D(\sqrt{H})=H^{1}$ and Theorem \ref{theorem1.2}, then we have
\begin{flalign}\label{2.3}
\|e^{itH}\phi\|_Z \leq C\|\phi\|_{H^1}.
\end{flalign}

The nonlinear Schr\"{o}dinger equation \eqref{equation1} with initial data at $t=s,u(s)=\phi_0$ can be written as the following integral equation
\begin{flalign}\label{2.4}
u(t)=e^{itH}\phi_0+\mathscr{P}_su(t).
\end{flalign}
We conclude that equation \eqref{2.4} has a unique solution in space $Z$ if $\|\phi_0\|_{H_1}$ is small enough because that $\mathscr{P}_s$ is a contraction on $Z(\delta)$. In particular, for $\phi_-\in H_{\delta}$ the unique solution that satisfies \eqref{1.2.2} is the fixed point when $s=-\infty$,
\begin{flalign}\label{2.5}
u(t)=e^{itH}\phi_-+\mathscr{P}_{-\infty}u(t)
\end{flalign}
and the unique $\phi_+$ that satisfies \eqref{1.2.3} is given by
\begin{flalign}\label{2.6}
\phi_+=\phi_-+\int_{-\infty}^{+\infty}e^{-itH}\arrowvert u\arrowvert^{p-1}u(t)dt.
\end{flalign}
Now our preparations are complete. Next, we can prove Theorem \ref{theorem1.3}.

\vspace{1cm}
\textbf{Proof of Theorem \ref{theorem1.3}:} 
\begin{proof}
By \eqref{2.5}, we obtain that for $\epsilon$ small enough
\begin{flalign}\label{2.7}
\|u\|_Z \leq \|e^{itH}\epsilon\phi_{-}\|_{Z}+c\|u\|_Z,
\end{flalign}
where there exist some $c<1$ since $\mathscr{P}_{-\infty}$ is contraction in $Z(\delta)$ for $\delta$ is small enough. Therefore, we have 
\begin{flalign}\label{2.8}
\|u\|_Z \leq C\epsilon\|e^{itH}\phi\|_Z.
\end{flalign}

Furthermore, from \eqref{assum1} in assumption (\uppercase\expandafter{\romannumeral1}) , estimates \eqref{2.1} and \eqref{2.6}, we have
\begin{flalign}\label{2.9}
\|(S_A-I)\epsilon\phi_-\|_{L^{1+p}}\leq C\int_{-\infty}^{+\infty}\frac{1}{|t|^{d}}\|u(t)\|^p_{L^{1+p}}dt
\end{flalign}
with $d=\dfrac{n(p-1)}{2(1+p)}<1$. Using \eqref{2.3} and \eqref{2.8}, for $|t|\leq 1$ we have
\begin{flalign*}
\int_{|t|\leq 1}\frac{1}{|t|^{d}}\|u(t)\|^p_{L^{p+1}}dt\leq C\|u(t)\|^p_{Z}\leq C\|e^{itH}\epsilon\phi_-\|^p_{Z}\leq \epsilon^p\|\phi_-\|_{H^1}^p,
\end{flalign*}
and for $|t|\geq 1$ we have
\begin{flalign*}
\int_{|t|\geq 1}\frac{1}{|t|^{d}}\|u(t)\|^p_{L^{p+1}}dt
\leq [\int_{|t|\geq 1}\frac{1}{|t|^{1+1/(r-p)}}dt]^{1-p/r}\|u(t)\|^p_{Z}\leq C\epsilon^p\|\phi_-\|_{H^1}^p.
\end{flalign*}
Thus, we obtain that
\begin{flalign}\label{2.10}
\|(S_A-I)\epsilon\phi_-\|_{L^{1+p}}\leq C\epsilon^p\|\phi_-\|_{H^1}^p,\ for\ \epsilon \ small\ enough.
\end{flalign}

On the other hand, we note that
\begin{flalign*}
|(|u|^{p-1}u-|w|^{p-1}w)|\leq C(|u|+|w|)^{p-1}|u-w|,\ for\ u,w\in C.
\end{flalign*}
Then the functional $u\rightarrow\int_{\mathbb{R}^{n}}|u|^{p+1}dx$ is bounded and uniformly Lipschitz continuous on bounded sets on $L^{1+p}$.
It follows from Theorem \ref{theorem1.2} that
\begin{flalign}\label{2.11}
E(u)=&\lim_{t\rightarrow\pm\infty}\int [\frac{1}{2}|\sqrt{H}u(t,.)|^2+|u(t,.)|^{p+1}]dx\nonumber\\
=&\frac{1}{2}\|\sqrt{H}\phi_{\pm}\|_{L^{2}}^2+\lim_{t\rightarrow \pm\infty}\int |e^{itH}\phi_\pm|^{p+1}dx,
\end{flalign}
where we used that $H^{1}\subset L^{1+p}$ and that since $D(\sqrt{H})=H^1$ the operator $(-\Delta+I)^{1/2}(H+I)^{-1/2}$ is bounded in $L^2$. Furthermore, from these facts we have
\begin{flalign*}
&\|e^{it_1H}\phi_\pm-e^{it_2H}\phi_\pm\|_{L^{1+p}}\\
\leq&C\|e^{it_1H}\phi_\pm-e^{it_2H}\phi_\pm\|_{H^1}\\
\leq&C\|(I-e^{i(t_2-t_1)H})\sqrt{H}\phi_\pm\|_{L^2}\rightarrow 0
\end{flalign*}
as $(t_2-t_1)\rightarrow 0$. Then $t\rightarrow \ e^{itH}\phi_\pm$ are uniformly continuous functions from $\mathbb{R}$ into $L^{1+p}$ and since also $e^{itH}\phi_\pm \in L^r(\mathbb{R},\,L^{1+p})$, we obtain that
\begin{flalign}\label{2.12}
\lim_{t\rightarrow\pm\infty}\|e^{itH}\phi_\pm\|_{L^{1+p}}=0.
\end{flalign}
It follows from \eqref{2.11} and \eqref{2.12} that
\begin{flalign*}
E(u)=\frac{1}{2}\|\sqrt{H}\phi_{\pm}\|_{L^{2}}^2
\end{flalign*}
and then there are constants $C_\pm >0$ such that
\begin{flalign}\label{2.13}
\|\phi_+\|_{H^1}\leq C_-\|\phi_-\|_{H^1}\leq C_+\|\phi_+\|_{H^1}.
\end{flalign}
Furthermore, as $\phi_+=S_A\phi_-$, we have
\begin{flalign}\label{2.14}
\frac{1}{\epsilon}\|(S_A-I)\epsilon\phi_-\|_{H^1}\leq C.
\end{flalign}

Note that $p>1$ it follows from \eqref{2.10} that $\dfrac{1}{\epsilon}(S_A-I)\epsilon\phi_-$ converges strongly to zero in $L^{1+p}$ and by \eqref{2.14} it is bounded in $H^1$. And since $L^{1+p}\cap H^1$ is dense in $H^1$, we have that $\frac{1}{\epsilon}(S_A-I)\epsilon\phi_-$ converges weakly to zero in $H^1$. Therefore for any $\phi,\psi\in H^{1}$, we obtain that
\begin{flalign}\label{2.16}
\lim_{\epsilon\downarrow 0}\frac{1}{\epsilon}[(S\epsilon\phi,\psi)-(S_L\epsilon\phi,\psi)]=\lim_{\epsilon\downarrow 0}\frac{1}{\epsilon}((S_A-I)W_-\epsilon\phi,W_+\psi)=0,
\end{flalign}
where we used $W_\pm$ are bounded operators from $H^1$ onto $H^1$. Therefore, the proof of Theorem \ref{theorem1.3} complete.
\end{proof}

Next we will prove Corollary \ref{corollary1.4} as follows. 
\begin{proof}
Let's denote $x=(x_{1},x_{2})\in\mathbb{R}^{2}$ with the vector $(x_{1},x_{2},0,0,...)$ in $\mathbb{R}^{n}$. For any $x\in\mathbb{R}^{n}$ and any $\phi_{0},\psi_{0}$ with $\hat{\phi}_{0},\hat{\psi}_{0}\in C_{0}^{\infty}$, we can define the following scattering states as in \cite{EW2}:
\begin{flalign}\label{1.4.1}
\phi_{\xi}:=e^{im\xi\cdot x}\phi_{0},
\end{flalign}
i.e.,
\begin{flalign}\label{1.4.2}
\hat{\phi}_{\xi}(m_{1})=\hat{\phi}_{0}(m_{1}-m\xi),
\end{flalign}
the $\hat{\psi}_{\xi}$ is similar defined. Thus it follows that 
\begin{flalign}\label{1.4.3}
|\xi|(i(S_{L}-I)\phi_{\xi},\,\psi_{\xi})=\int_{-\infty}^{+\infty}(A(x+\tau\hat{\xi})\phi_{0},\,\psi_{0})d\tau+O(|\xi|^{-1}),\ as\, |\xi|\rightarrow\infty.
\end{flalign}
Since the set of $\phi_{0},\psi_{0}$ with $\hat{\phi}_{0},\hat{\psi}_{0}\in C_{0}^{\infty}$ is dense in $L^{2}$ one reconstructs $A$ uniquely as an operator and as a function for every $x$ because $A(x)$ is continuous. By \eqref{1.4.3} from the high velocity limit of $S_{L}$ one reconstructs the Radon transform of the matrix elements of $A$ in a dense set of states. Inverting this Radon transform one uniquely reconstructs $A$. For further details see \cite{EW2}. Therefore, the proof of Corollary \ref{corollary1.4} is completed.
\end{proof}

\section{Proof of Theorem \ref{theorem1.5}}
In this section, we consider Schr\"odinger equation \eqref{equation1} with $\phi \in \Sigma$ in $n=1$. Let potentials $A(x)$ and $V(x)$ satisfy assumption (\uppercase\expandafter{\romannumeral3}). It is worth noting that we don't need the assumptions (\uppercase\expandafter{\romannumeral1}) and (\uppercase\expandafter{\romannumeral2}). We use the decay estimate in \cite{WD} instead of using the dispersive estimate \eqref{assum1}, the condition to potential $A(x)$ with \eqref{assum3} instead of \eqref{assum2}. We will use the same notation to distinguish the difference from the proof of Theorem \eqref{theorem1.3}.

Let's denote $Z=L^r(\mathbb{R},\,X)\cap L^{\infty}(\mathbb{R},\,X)$, where $r=\frac{p-1}{1-n}$. The norm of $Z$ is given by
\begin{flalign*}
\|u\|_Z=\|u\|_{L^r(\mathbb{R},\, X)}+\|u\|_{L^{\infty}(\mathbb{R},\,X)}.
\end{flalign*}
Applying the decay estimate in our other paper \cite{WD}, we have 
\begin{flalign}\label{3.1}
\|e^{itH}\phi\|_{L^\infty(\mathbb{R}^{n})}\leq C|t|^{-n/2}\|\phi\|_{\Sigma(\mathbb{R}^{n})},
\end{flalign}
where $\Sigma=\{\phi\in\mathscr{S'}(\mathbb{R}^{n}),\ (-\Delta)^{\frac{s}{2}}\phi\in L^{2}(\mathbb{R}^{n}),\ \arrowvert x\arrowvert^{s}\phi\in L^{2}(\mathbb{R}^{n}),\ |x|^{s_{2}}(-\Delta)^{\frac{s_{1}}{2}}\phi\in L^{2}(\mathbb{R}^{n})\}$. Moreover, note that the global existence of the solution of Schr\"odinger equation \eqref{equation1} in \cite{B} as follows  
\begin{flalign}\label{3.2}
\|e^{itH}\phi\|_{\Sigma}\leq C.
\end{flalign}
Thus, combining \eqref{3.1} and \eqref{3.2}, denoting $X=L^{\infty}\cap\Sigma$, we obtain that 
\begin{flalign}\label{3.3}
\|e^{itH}f\|_{X}\leq C\langle t\rangle^{-n\theta/2}\|f\|_{X},\ 0\leqslant\theta\leqslant1.
\end{flalign}
Meanwhile, we estimate that
\begin{flalign*}
&\||u|^{p-1}u\|_{\Sigma}\\
=&(\||u|^{p-1}u\|^{2}_{\dot{H}^{s}}+\||x|^{s}|u|^{p-1}u\|^{2}_{L^{2}}+\||x|^{s_{2}}(-\Delta)^{\frac{s_{1}}{2}}(|u|^{p-1}u)\|^{2}_{L^{2}})^{\frac{1}{2}}\\
\leqslant&(C\||u|^{2(p-1)}_{L^{\infty}}\|u\|^{2}_{\dot{H}^{s}}+C\||u|^{2(p-1)}_{L^{\infty}}\||x|^{s}u\|^{2}_{L^{2}}+C\||u|^{2(p-1)}_{L^{\infty}}\||x|^{s_{2}}(-\Delta)^{\frac{s_{1}}{2}}u\|^{2}_{L^{2}})^{\frac{1}{2}}\\
\leqslant&C\||u|^{p-1}_{L^{\infty}}\cdot\|u\|_{\Sigma}\\
\leqslant&C\||u|^{p}_{X}.
\end{flalign*}

We claim that $\mathscr{P}_s$ is a contraction on $Z(\delta_{1})=\{u\in Z, \|u\|_Z\leq \delta_{1}\}$ if $\delta_{1}\ll1$. In fact, on the one hand
\begin{flalign*}
\|\mathscr{P}_su-\mathscr{P}_sv\|_{X}=&\|\int_{s}^{t}e^{i(t-\tau)H}(|u|^{p-1}u-|v|^{p-1}v)d\tau\|_{X}\\
\leqslant&C\int_{s}^{t}\|e^{i(t-\tau)H}(|u|^{p-1}u-|v|^{p-1}v)\|_{X}d\tau\\
\leqslant&C\int_{s}^{t}\langle t-\tau\rangle^{-\frac{n\theta}{2}}\||u|^{p-1}u-|v|^{p-1}v\|_{\Sigma}d\tau\\
\leqslant&C\int_{s}^{t}\langle t-\tau\rangle^{-\frac{n\theta}{2}}(\|u\|_{X}+\|v\|_{X})^{p-1}\cdot\|u-v\|_{X}d\tau\\
\leqslant&C(\|u\|_{Z}+\|v\|_{Z})^{p-1}\cdot\|u-v\|_{Z}.
\end{flalign*}
The last inequality holds since that choose $I_{1}=[s,t]\bigcap[t-1,t+1]$, $I_{2}=[s,t]\setminus[t-1,t+1]$, we have
\begin{flalign*}
\|\mathscr{P}_su\|_{X}\leqslant& C\int_{s}^{t}\langle t-\tau\rangle^{-\frac{n\theta}{2}}\|u\|^{p}_{X}d\tau\\
\leqslant&C\int_{I_{1}}\langle t-\tau\rangle^{-\frac{n\theta}{2}}d\tau\cdot\|u\|^{p}_{L_{t}^{\infty}(\mathbb{R},\,X)}\\
&+C(\int_{I_{2}}\langle t-\tau\rangle^{1+\frac{1}{r-p}}d\tau)^{\frac{1}{r}}\cdot(\int_{I_{2}}\|u\|_{X}^{p\cdot\frac{r}{p}})^{\frac{p}{r}}\\
\leqslant&C\|u\|^{p}_{Z}.
\end{flalign*}
On the other hand, we have 
\begin{flalign*}
\int_{-\infty}^{\infty}\|\mathscr{P}_su-\mathscr{P}_sv\|_{X}^{r}dt\leqslant&C\int_{\mathbb{R}}\arrowvert\int_{s}^{t}\langle t-\tau\rangle^{-\frac{n\theta}{2}}(\|u\|_{X}+\|v\|_{X})^{p-1}\cdot\|u-v\|_{X}d\tau\arrowvert^{r}dt\\
=&C\int_{\mathbb{R}}\arrowvert\langle t-\tau\rangle^{-\frac{n\theta}{2}}\ast((\|u\|_{X}+\|v\|_{X})^{p-1}\|u-v\|_{X})\arrowvert^{r}dt\\
\leqslant&C(\int_{\mathbb{R}}\arrowvert(\|u\|_{X}+\|v\|_{X})^{p-1}\cdot\|u-v\|_{X}\arrowvert^{\frac{r}{p}}dt)^{\frac{p}{r}}\\
\leqslant&C(\|u\|_{Z}+\|v\|_{Z})^{p-1}\|u-v\|_{Z}.
\end{flalign*}
Combining above, we obtain
\begin{flalign}\label{3.4}
\|\mathscr{P}_su-\mathscr{P}_sv\|_{Z}\leqslant C(\|u\|_{Z}+\|v\|_{Z})^{p-1}\|u-v\|_{Z}.
\end{flalign}
Therefore, we complete the claim.

The nonlinear Schr\"{o}dinger equation with initial data at $t=s,u(s)=\phi_0$ can be written as the following integral equation
\begin{flalign}\label{3.5}
u(t)=e^{itH}\phi_0+\mathscr{P}_su(t).
\end{flalign}
Since $\mathscr{P}_s$ is a contraction on $Z(\delta_{1})$, \eqref{3.5} has a unique solution in $Z$ if $||\phi_0||_{\Sigma}$ is small enough.
In particular for $\phi_-\in\Sigma$ the unique solution that satisfies \eqref{1.2.2} is the fixed point when $s=-\infty$,
\begin{flalign}\label{3.6}
u(t)=e^{itH}\phi_-+\mathscr{P}_{-\infty}u(t)
\end{flalign}
and the unique $\phi_+$ that satisfies \eqref{1.2.3} is given by
\begin{flalign}\label{3.7}
\phi_+=\phi_-+\int_{-\infty}^{+\infty}e^{-itH}|u|^{p-1}u(t)dt.
\end{flalign}
By \eqref{3.2} for $\epsilon$ small enough, we have
\begin{flalign}\label{3.8}
\|u\|_Z \leq \|e^{itH}\epsilon\phi\|_Z+c\|u\|_Z,
\end{flalign}
where there exsit some $c<1$ since $\mathscr{P}_{-\infty}$ is contraction in $Z(\delta_{1})$ for $\delta_{1}$ is small enough. Thus from \eqref{3.8},
\begin{flalign}\label{3.9}
\|u\|_Z \leq C\epsilon\|e^{itH}\phi\|_Z.
\end{flalign}

Using \eqref{3.1} and \eqref{3.7}, we have
\begin{flalign}\label{3.10}
\|(S-I)\epsilon\phi_-\|_{L^{\infty}} &\leq C\int_{-\infty}^{+\infty}\frac{1}{|t|^{\frac{n}{2}}}\||u|^{p-1}u(t)\|_{\Sigma}dt\nonumber\\
&\leq C\int_{-\infty}^{+\infty}\frac{1}{|t|^{\frac{n}{2}}}\|u(t)\|^{p}_{X}dt.
\end{flalign}
From \eqref{3.3} and \eqref{3.9}, for $|t|\leq 1$ we have
\begin{flalign*}
\int_{|t|\leq 1}\frac{1}{|t|^{\frac{n}{2}}}\|u(t)\|^{p}_{X}dt\leq C\|u(t)\|^{p}_{Z}\leq C\|e^{itH}\epsilon\phi_-\|^{p}_{Z}\leq C\epsilon^{p}\|\phi_-\|_{X}^p.
\end{flalign*}
For $|t|\geq 1$ we have
\begin{flalign*}
\int_{|t|\geq 1}\frac{1}{|t|^{\frac{n}{2}}}\|u(t)\|^{p}_{X}dt
\leq C(\int_{|t|\geq 1}|t|^{1+1/(r-p)}dt)^{1/r}\cdot\|u(t)\|^p_{Z}\leq C\epsilon^p\|\phi_-\|_{X}^p.
\end{flalign*}
Thus we have
\begin{flalign}\label{3.11}
\|(S-I)\epsilon\phi_-\|_{L^{\infty}}\leq C\epsilon^p\|\phi_-\|_{X}^p,\ for\ \epsilon \ small\ enough.
\end{flalign}

Furthermore, 
\begin{flalign}\label{3.12}
|(|u|^{p-1}u-|w|^{p-1}w)|\leq C(|u|+|w|)^{p-1}|u-w|,\ for\ u,w\in C.
\end{flalign}
Then the functional $u\rightarrow\int_{\mathbb{R}} |\phi|^{p+1}dx$ is bounded and uniformly Lipschitz continuous on bounded sets on $L^{1+p}$. It follows from Theorem \ref{theorem1.2} that
\begin{flalign}\label{3.13}
E(u)=&\lim_{t\rightarrow\pm\infty}\int [\frac{1}{2}|\sqrt{H}u(t,.)|^2+|u(t,.)|^{p+1}]dx\nonumber\\
=&\frac{1}{2}\|\sqrt{H}\phi_{\pm}\|_{L^{2}}^2+\lim_{t\rightarrow \pm\infty}\int |e^{itH}\phi_\pm|^{p+1}dx,
\end{flalign}
where we used that since $D(\sqrt{H})=H^1$ the operator $(-\Delta+I)^{1/2}(H+I)^{-1/2}$ is bounded in $L^2$. Meanwhile,
\begin{flalign*}
&\|e^{it_1H}\phi_\pm-e^{it_2H}\phi_\pm\|_{\Sigma}\\
\leq&C\|e^{it_1H}\phi_\pm-e^{it_2H}\phi_\pm\|_{H^1}\\
\leq&C\|(I-e^{i(t_2-t_1)H})\sqrt{H}\phi_\pm\|_{L^2}\rightarrow 0
\end{flalign*}
as $(t_2-t_1)\rightarrow 0$. Then $t\rightarrow \ e^{itH}\phi_\pm$ are uniformly continuous functions from $\mathbb{R}$ into $\Sigma$ and since also $e^{itH}\phi_\pm \in L^r(\mathbb{R},\,\Sigma)$, we have that
\begin{flalign}\label{3.13-1}
\lim_{t\rightarrow\pm\infty}||e^{itH}\phi_\pm||_{\Sigma}=0.
\end{flalign}
It follows from \eqref{3.13} and \eqref{3.13-1} that
\begin{flalign*}
E(u)=\frac{1}{2}\|\sqrt{H}\phi_{\pm}\|_{L^{2}}^2,
\end{flalign*}
and then there are constants $C_\pm >0$ such that
\begin{flalign}\label{3.14}
\|\phi_+\|_{H^1}\leq C_-\|\phi_-\|_{H^1}\leq C_+\|\phi_+\|_{H^1}.
\end{flalign}
Therefore, as $\phi_+=S_A\phi_-$ and the definition of $\Sigma$, we have
\begin{flalign}\label{3.15}
\frac{1}{\epsilon}||(S-I)\epsilon\phi_-||_{\Sigma}\leq C.
\end{flalign}

Since $p>1+\frac{2}{n}$ it follows from \eqref{3.11} that $\dfrac{1}{\epsilon}(S-I)\epsilon\phi_-$ converges strongly to zero in $\Sigma$ and by \eqref{3.15} it is bounded in $\Sigma$. As $L^{\infty}\cap \Sigma$ is dense in $\Sigma$ we have that $\frac{1}{\epsilon}(S-I)\epsilon\phi_-$ converges weakly to zero in $\Sigma$. Then for any $\phi,\psi\in \Sigma$
\begin{flalign*}
\lim_{\epsilon\downarrow 0}\frac{1}{\epsilon}[(S\epsilon\phi,\psi)-(S_L\epsilon\phi,\psi)]=
\lim_{\epsilon\downarrow 0}\frac{1}{\epsilon}((S-I)W_-\epsilon\phi,W_+\psi)=0,
\end{flalign*}
where we used that $W_\pm$ are bounded operators from $\Sigma$ onto $\Sigma$.\\

Next we will prove Corollary \ref{corollary1.6} as follows. 
\begin{proof}
Similar the proof of Corollary \ref{corollary1.4}, we can invert the Radon transform with megnetic potentials one reconstructs $A(x)$ and $V(x)$. Since the set of $\phi_{0},\psi_{0}$ with $\hat{\phi}_{0},\hat{\psi}_{0}\in C_{0}^{\infty}$ is dense in $L^{2}$ one reconstructs $A(x)$ and $V(x)$ uniquely as operators and as functions for every $x$ because $A(x)$ and $V(x)$ are continuous. This implies Corollary \ref{corollary1.6} since $A$ and $V$ are continuous functions of $x$ it is possible to reconstruct $A(x)$ and $V(x)$ directly as functions from \eqref{1.4.3}. Therefore, the proof of Corollary \ref{corollary1.6} is completed.
\end{proof}

\end{document}